\documentclass[runningheads]{lncse}

\usepackage{graphicx}
\usepackage{amsmath,amssymb}
\usepackage{amsfonts}
\usepackage{algorithm}
\usepackage{algorithmic}
\usepackage{subfigure}
\graphicspath{ {./eps/} }

\begin{document}
%\linenumbers
\title{Time-Space Adaptive Method of Time Layers for the Advective  Allen-Cahn Equation}

\titlerunning{Advective Allen-Cahn equation}

\author{ Murat Uzunca\inst{1}  \and B\"{u}lent Karas\"{o}zen\inst{2} \and Ay\c{s}e Sar{\i}ayd{\i}n Filibelio\u{g}lu  \inst{1}}

\authorrunning{Uzunca et al.}

\institute{
Institute of Applied Mathematics, Middle East Technical University, 06800 Ankara, Turkey,  {\tt uzunca@gmail.com},  {\tt saayse@metu.edu.tr} \and
Department of Mathematics \& Institute of Applied Mathematics, Middle East Technical University, 06800 Ankara, Turkey, {\tt bulent@metu.edu.tr}
}

\maketitle

\begin{abstract}

We develop an adaptive method of time layers with a linearly implicit Rosenbrock method as time integrator and symmetric interior penalty Galerkin method for space discretization for the advective Allen-Cahn equation with non-divergence-free velocity fields. Numerical simulations for convection dominated problems demonstrate the accuracy and efficiency of the adaptive algorithm for  resolving the sharp layers occurring in interface problems with small surface tension.
\end{abstract}

%%%%%%%%%%%%%%%%%%%%%%%%%%%%%%%%%%%%%%%%%%%%%%%%%%%%%%%%%%%%%%%%%%%%%%%%%%%%%%%
%%%%%%%%%%%%%%%%%%%%%%%%%%%%%%%%%%%%%%%%%%%%%%%%%%%%%%%%%%%%%%%%%%%%%%%%%%%%%%%
\section{Introduction}

Interfacial dynamics has great importance in  modeling of multi phase flow in material sciences, and binary fluids flow movement. We consider the Allen-Cahn equation with advection as a model of diffuse interface for two phase flows \cite{liu12}
\begin{equation} \label{advac}
 \frac{\partial u}{\partial t} = \mathcal{L} u - \frac{1}{\epsilon}f(u)\quad  \hbox{in } \Omega \times (0,T]
\end{equation}
under homogeneous Neumann boundary conditions, where $\mathcal{L}$ denotes the linear operator related to the diffusion and advection parts of the system, i.e. $\mathcal{L}u=\epsilon\Delta u - \nabla\cdot (\mathbf{V}u)$. The term $f(u)=F'(u)=2u(1 - u)(1 -2u)$ characterizes the cubic bistable nonlinearity with the double--well potential $F(u)$ of the two phases, and $\epsilon$ describes the surface tension. We consider a prescribed fixed velocity field $\mathbf{V}=(V_1, V_2)^T$. In coupled incompressible fluid mechanics and  diffusive interface models, the velocity field satisfies the Navier--Stokes equations \cite{liu12}, and therefore is  divergence free, i.e. $\nabla \cdot \mathbf{V} = 0$. We consider in this work the non-divergence-free velocity fields which are either expanding ($\nabla \cdot \mathbf{V}> 0$) or sheering ($\nabla \cdot \mathbf{V} < 0$) as in \cite{liu12}.

The advective Allen-Cahn equation \eqref{advac} describes the diffuse interface dynamics associated with surface energies, and has two different time scales; the small surface tension, and the convection time scale. Both time scales cause computational stiffness \cite{liu12}. The dynamics of surface tension in two-phase fluids are studied numerically by  the level--set algorithm method and the diffuse interface method \cite{liu12}.

In this work we apply the adaptive method of time layers (AMOT) \cite{deuflhard12}, or adaptive Rothe method, where the advective Allen-Cahn equation \eqref{advac} is discretized   first in time then in space, in contrast to the usual method of lines approach. Hereby spatial discretization is considered as a perturbation of the time integration. AMOT was applied to linear and nonlinear partial differential equations using linearly implicit time integrators in several papers \cite{bornemann90,deuflhard12,frank11,lang93,lang01,lang01a}. We have chosen the linearly three stage Rosenbrock (ROS3P) method \cite{lang01} as the time integrator. ROS3P solver is third order in time, L-stable and can efficiently deal with the large stiff systems arising from the discretization of \eqref{advac}. It does not show any order reduction in time in contrast to other Rosenbrock methods of order higher than two \cite{lang01}. Unlike the fully implicit schemes, it requires only the solution of  three linear systems per time step with the same coefficient matrix. In non-stationary models, the potential internal/boundary layers moves as the time progresses. The time step-sizes have to be adapted properly to resolve these layers accurately. The simple  embedded  a posteriori error estimator as the difference of second and third order ROS3P solvers allow the construction of   an efficient adaptive time integrator. To resolve the sharp layers and oscillations in advection-dominated regimes, we apply symmetric interior penalty Galerkin (SIPG) method \cite{arnold01,riviere08}, as a stable space discretization in the family of discontinuous Galerkin (dG) methods. Further, we apply the adaptive SIPG method in space with the residual-based a posteriori  error estimator \cite{schotzau09,uzunca14adg} to handle unphysical oscillations. The spatial mesh is refined or coarsened  locally to obtain an accurate approximation with less degree of freedoms (DoFs) and less computational time. We show in numerical experiments that the proposed time-space algorithm AMOT is capable of damping  the oscillations which may vary as the time progresses.

The paper is organized as follows. In Section~\ref{dg} we give the fully discrete formulation of the advective AC model \eqref{advac}. The time-space adaptive algorithm is described in Section~\ref{adap}. In Section~\ref{numeric}, results of numerical experiments for convection dominated expanding and sheering flows are presented.

%%%%%%%%%%%%%%%%%%%%%%%%%%%%%%%%%%%%%%%%%%%%%%%%%%%%%%%%%%%%%%%%%%%%%%%%%%%%%%%
%%%%%%%%%%%%%%%%%%%%%%%%%%%%%%%%%%%%%%%%%%%%%%%%%%%%%%%%%%%%%%%%%%%%%%%%%%%%%%%
\section{Time-Space Discretization}
\label{dg}

In this section we apply the method of time layers to discretize the model \eqref{advac} in time. The resulting sequence of elliptic problems are discretized by the SIPG method at each time step. We consider the partition of time interval $[0,T]$, as $I_k=(t^{k-1},t^k]$ with  time step-sizes $\tau_k=t^{k}-t^{k-1}$, $k=1,2, \ldots, J$. The approximate solution at the time $t=t^k$ is denoted by $u^k\approx u(t^k)$. We apply the 3-stage Rosenbrock solver ROS3P \cite{lang01} with an embedded error estimator in time:
\begin{equation}\label{ros_system}
\left(\frac{1}{\gamma\tau_k}- J^{k-1}\right)K_i = \mathcal{L}\left( z^i\right) - f\left(z^i\right) + \sum_{j=1}^{i-1} \frac{c_{ij}}{\tau_k}K_j, \quad j=1,2,3,
\end{equation}
where $z^i=u^{k-1} +\sum_{j=1}^{i-1} a_{ij}K_j$ and $J^{k-1}:=J(u^{k-1})$ is the Jacobian $J(u)=\partial_u({\mathcal L}u-f(u))$
at $u^{k-1}$. The second order solution $\hat{u}^k$ and the third order solution $u^k$ are given by
\begin{equation}\label{sol}
\begin{aligned}
\hat{u}^k &=  u^{k-1} + \hat{m}_1K_1+\hat{m}_2K_2+\hat{m}_3K_3 ,\\
u^k &=  u^{k-1} + m_1K_1+m_2K_2+m_3K_3
\end{aligned}
\end{equation}
with the same stage vectors $K_i$. For the derivation of ROS3P solver and for parameter values, we refer to  \cite{lang01}.
The difference of  the solutions $u^k$  and $\hat{u}^k$ is used as an error indicator in the time-adaptivity. Due to the linearly  implicit nature of the Rosenbrock methods, the stage vectors $K_i$ in \eqref{ros_system} are solved using linear systems with the same coefficient matrix, which increases the computational efficiency in time integration of nonlinear PDEs \cite{deuflhard12,lang01a}.

The semi-discrete systems \eqref{ros_system} are discretized in space by the SIPG method with upwinding for the convective term \cite{ayuso09}. On the time interval  $I_n=(t^{k-1},t^k]$, we  consider a family ${\mathcal{T}}_h^k$ of shape regular elements (triangles) $E\in{\mathcal{T}}_h^k$, and we denote the initial one by ${\mathcal{T}}_h^0$. Here, the mesh ${\mathcal{T}}_h^k$ is obtained by local refinement/coarsening of the mesh ${\mathcal{T}}_h^{k-1}$ from the previous time step. Then, with the dG finite element space $V_h^k:=V_h(\mathcal{T}_h^k)$, on an individual time step $I_n=(t^{k-1},t^k]$, in time ROS3P and in space SIPG discretized fully discrete scheme for the model \eqref{advac} reads as: for all $v_h^k\in V_h^k$, find $u_h^k$ (or $\hat{u}_h^k$) in \eqref{sol} with $K_i\in V_h^k$, $i=1,2,3$, satisfying
\begin{equation}\label{fully}
\left( \left(\frac{1}{\gamma\tau_k}- J_h^{k-1}\right)K_i,v_h^k\right) =  -a_h(z_h^i,v_h^k) - b_{h}(z_h^i, v_h^k) + \left( \sum_{j=1}^{i-1} \frac{c_{ij}}{\tau_k}K_j,v_h^k \right),
\end{equation}
where $z_h^i=u_h^{k-1} +\sum_{j=1}^{i-1} a_{ij}K_j$, $J_h^{k-1}=J(u_h^{k-1})$ and $(\cdot , \cdot)$ stands for the discrete inner product $(\cdot , \cdot)_{L^2(\mathcal{T}_h^k)}$. The forms $a_h(u_h^k,v_h^k)$ and $b_h(u_h^k,v_h^k)$ are the bilinear and linear forms given by
\begin{align*}
a_{h}(u_h^k, v_h^k)=& \sum \limits_{E \in {\mathcal{T}}_h^k} \int_{E} \epsilon \nabla u_h^k\cdot\nabla v_h^k dx + \sum \limits_{E \in {\mathcal{T}}_h^k} \int_{E} (\mathbf{V}\cdot \nabla u_h^k + (\nabla\cdot\mathbf{V}) u_h^k) v_h^k dx \\
&+ \sum \limits_{E \in {\mathcal{T}}_h^k}\int_{\partial E^-\setminus\Gamma_h^- } \mathbf{V}\cdot \mathbf{n}_E ((u_{h}^{out})^k-(u_{h}^{in})^k)  v_h^k ds \\
&- \sum \limits_{E \in {\mathcal{T}}_h^k} \int_{\partial E^-\cap \Gamma_h^{-}} \mathbf{V}\cdot \mathbf{n}_E (u_{h}^{in})^k v_h^k ds + \sum \limits_{ e \in \Gamma_h^k}\frac{\sigma \epsilon}{h_{e}} \int_{e} [u_h^k]\cdot[v_h^k] ds  \\
&- \sum \limits_{ e \in \Gamma_h^k} \int_{e} ( \{\epsilon \nabla v_h^k \}\cdot[u_h^k] + \{\epsilon \nabla u_h^k \}\cdot [v_h^k] )ds,\\
b_{h}(u_h^k, v_h^k) =& \sum \limits_{E \in {\mathcal{T}}_h^k} \int_{K} \frac{1}{\epsilon}f(u_h^k) v_h^k dx,
\end{align*}
where $u_{h}^{out}$ and $u_{h}^{in}$ denote the traces on an edge from outside and inside of an element $E$, respectively, $h_e$ is the length of an edge $e$, $\Gamma_h^k$ is the set of interior edges, $\partial E^-$ and $\Gamma_h^-$ are the sets of inflow boundary edges of an element $E\in{\mathcal{T}}_h^k$ and on the boundary $\partial\Omega$, respectively. The parameter $\sigma$ is called the penalty parameter to penalize the jumps in dG schemes, and $[\cdot]$ and $\{\cdot\}$ stand as the jump and average operators, respectively \cite{riviere08}.

%%%%%%%%%%%%%%%%%%%%%%%%%%%%%%%%%%%%%%%%%%%%%%%%%%%%%%%%%%%%%%%%%%%%%%%%%%
%%%%%%%%%%%%%%%%%%%%%%%%%%%%%%%%%%%%%%%%%%%%%%%%%%%%%%%%%%%%%%%%%%%%%%%%%%%%%%%
\section{Adaptive Method of Time Layers (AMOT)}
\label{adap}

The goal of AMOT is on each time step $I_k=(t^{k-1},t^k]$ adjusting the time step-size and the spatial mesh adaptively. To do this, AMOT aims to bound the total error $\|u(t^k)-\hat{u}_h^k\|$ by suitable temporal and spatial estimators, where $u(t^k)$ is the true solution of the continuous model \eqref{advac} and $\hat{u}_h^k$ is the $2^{nd}$ order (in time) discrete solution of the fully discrete system \eqref{fully} on ${\mathcal{T}}_h^{k-1}$, at the time $t=t^k$. In order to define the temporal and spatial estimators separately, we replace the true solution $u(t^k)$ by its best available approximation $\overline{u_h^{k,+}}$ which is the $3^{rd}$ order (in time) discrete solution of the fully discrete system \eqref{fully} on an auxiliary very fine mesh $\overline{{\mathcal{T}}_h^k}\supset{\mathcal{T}}_h^{k-1}$, and we add and subtract  the term $u_h^k$ which is the $3^{rd}$ order (in time) discrete solution of the fully discrete system \eqref{fully} on ${\mathcal{T}}_h^{k-1}$ at the time $t=t^k$. Then, similar to \cite[Sec. 9.2]{deuflhard12}, we get
\begin{equation}\label{error}
\begin{aligned}
\|u(t^k)-\hat{u}_h^k\|_{L^2({\mathcal{T}}_{h}^{k-1})} & \approx \|\overline{u_h^{k,+}}-\hat{u}_h^k\|_{L^2({\mathcal{T}}_{h}^{k-1})} = \|\overline{u_h^{k,+}}-u_h^k+u_h^k-\hat{u}_h^k\|_{L^2({\mathcal{T}}_{h}^{k-1})} \\
& \leq \underbrace{\|\overline{u_h^{k,+}}-u_h^k\|_{L^2({\mathcal{T}}_{h}^{k-1})}}_{:=\varepsilon_S} + \underbrace{\|u_h^k-\hat{u}_h^k\|_{L^2({\mathcal{T}}_{h}^{k-1})}}_{:=\varepsilon_T}\\
& \leq TOL_S + TOL_T \leq TOL
\end{aligned}
\end{equation}
for a user prescribed tolerance $TOL$, and further, we set $TOL_T=\alpha TOL$ and $TOL_S=(1-\alpha) TOL$ for user defined $0<\alpha <1$. In \eqref{error}, the term $\varepsilon_T$ controls the temporal adjustment, while the term $\varepsilon_S$ controls the acceptance of spatial mesh. Note that the temporal error estimator $\varepsilon_T$ is nothing but the difference of the order 2 and order 3 (in time) solutions of the fully discrete system \eqref{fully} on ${\mathcal{T}}_h^{k-1}$ at the time $t=t^k$. As a result, on each time step $I_k$, AMOT starts on the spatial mesh ${\mathcal{T}}_h^{k-1}$ with the adjustment of the time step-size $\tau_k$ according to the acceptance relation \cite{deuflhard12}
\begin{equation}\label{tau}
\tau^* = \sqrt[3]{\frac{\rho TOL_T}{\varepsilon_T}}\tau_k
\end{equation}
with a safety factor $\rho\approx 0.9$, and the computed time step-size $\tau^*$ is accepted if $\varepsilon_T\leq TOL_T$.

\begin{algorithm}
\caption{AMOT Algorithm on a single time step $I_k=(t^{k-1},t^k]$}
\textbf{Input:} $u_h^{k-1}$, $\tau^*$, ${\mathcal{T}}_h^{k-1}$, $TOL_S$, $TOL_T$ \\
\textbf{Output:} $u_h^{k,+}$, $\tau_{k}$, $\tau^*$, ${\mathcal{T}}_h^k$
\begin{algorithmic}
\STATE {\bf do}\\
		  \qquad $\tau_k=\tau^*$\\
			\qquad compute $u_h^k$ and $\hat{u}_h^k$ on ${\mathcal{T}}_h^{k-1}$\\
			\qquad {\bf if} $\varepsilon_T>TOL_T$\\
			\qquad \qquad compute new step-size $\tau^*$ according to \eqref{tau}\\
			\qquad {\bf end if}\\
			\qquad compute error indicator $\eta$ and construct the auxiliary fine mesh $\overline{{\mathcal{T}}_h^k}$\\
			\qquad compute the best available approximation $\overline{u_h^{k,+}}$ on $\overline{{\mathcal{T}}_h^k}$\\
			\qquad {\bf if} $\varepsilon_S>TOL_S$\\
			\qquad \qquad refine elements $E\in {\mathcal{T}}_h^{k-1}$ with $(\varepsilon_S)_E>0.005\times TOL_S$\\
			\qquad \qquad coarsen elements $E\in {\mathcal{T}}_h^{k-1}$ with $(\varepsilon_S)_E<10^{-13}$\\
			\qquad \qquad construct the new spatial mesh ${\mathcal{T}}_h^k$\\
			\qquad {\bf end if}		\\
			\qquad compute $u_h^{k,+}$ on ${\mathcal{T}}_h^k$		
\STATE  {\bf until} $\varepsilon_T\leq TOL_T$ and $\varepsilon_S\leq TOL_S$
\end{algorithmic}
\label{algorithm}
\end{algorithm}

After time step-size adjustment, AMOT continues with the refinement and coarsening of the spatial mesh ${\mathcal{T}}_h^{k-1}$ to obtain the new spatial mesh ${\mathcal{T}}_h^k$ according to the spatial estimator $\varepsilon_S=\sum_{E} (\varepsilon_S)_E$ in \eqref{error}, where the local elements $E\in{\mathcal{T}}_h^{k-1}$ are refined for large $(\varepsilon_S)_E$ and the ones are coarsened for small $(\varepsilon_S)_E$. To match the elements $E\in{\mathcal{T}}_h^{k-1}$ to be refined, we check the condition $(\varepsilon_S)_E>0.005\times TOL_S$, while we check the condition $(\varepsilon_S)_E<10^{-13}$ to match the elements $E\in{\mathcal{T}}_h^{k-1}$ to be coarsened. Yet, to compute the spatial estimator $\varepsilon_S$, we need the best available approximation $\overline{u_h^{k,+}}$ which is the solution of the discrete system \eqref{fully} on a very fine auxiliary mesh $\overline{{\mathcal{T}}_h^k}\supset{\mathcal{T}}_h^{k-1}$. To construct the auxiliary fine mesh $\overline{{\mathcal{T}}_h^k}$, one needs a local error indicator to match the elements $E\in{\mathcal{T}}_h^{k-1}$ to be refined. We use residual-based error indicator \cite{schotzau09}
\begin{equation}\label{ind}
\eta =\left( \sum \limits_{E \in {\mathcal{T}}_h^{k-1}}\eta_E^2\right)^{1/2} \; , \quad \eta_E^2= \eta_{E_R}^2  + \eta_{E_0}^2 + \eta_{E_{\partial}}^2,
\end{equation}
where $\eta_{E_R}$ denote the cell residuals
\begin{equation*}
\eta_{E_R}^2 = \lambda_E^2\left\| \frac{u_h^k-u_h^{k-1}}{\tau_k} - \epsilon\Delta u_h^k + \nabla\cdot (\mathbf{V}u_h^k)  + \frac{1}{\epsilon}f(u_h^k)\right\|_{L^2(E)}^2,
\end{equation*}
for a weight function $\lambda_E$, while $\eta_{E_0}$ and $\eta_{E_{\partial}}$ stand for the edge residuals coming from the jump of the numerical solution on the interior and Neumann boundary edges, respectively, see \cite{schotzau09,uzunca14adg} for details. Using the local error indicators $\eta_E$ in \eqref{ind}, we construct the auxiliary fine mesh $\overline{{\mathcal{T}}_h^k}$ by refining the elements $E\in M_E\subset{\mathcal{T}}_h^{k-1}$, where the set $M_E$ is determined by the bulk criterion
$$
\sum \limits_{E \in M_E}\eta_E^2 \geq  \theta\sum \limits_{E \in {\mathcal{T}}_h^{k-1}}\eta_E^2.
$$
for a user prescribed $0<\theta <1$, where larger $\theta$ leads to more elements to be refined. In our simulations we take $\theta =0.9$ since we need a very fine auxiliary mesh. The AMOT procedure, see Algorithm~\ref{algorithm}, continues until the temporal and spatial acceptance conditions $\varepsilon_S\leq TOL_S$ and $\varepsilon_T\leq TOL_T$ are satisfied.

%%%%%%%%%%%%%%%%%%%%%%%%%%%%%%%%%%%%%%%%%%%%%%%%%%%%%%%%%%%%%%%%%%%%%%%%%%%%%%%
%%%%%%%%%%%%%%%%%%%%%%%%%%%%%%%%%%%%%%%%%%%%%%%%%%%%%%%%%%%%%%%%%%%%%%%%%%%%%%%
\section{Numerical Experiments}
\label{numeric}

In this section, we demonstrate the accuracy and efficiency of the proposed AMOT for expanding and sheering flow examples. In all examples, we set the tolerance $TOL=0.001$, the parameter $\alpha =0.5$ and the diffusion coefficient $\epsilon =0.01$. The spatial domain is taken as $\Omega =[-1,1]^2$ and the time interval is $[0,06]$. For the SIPG discretization we use piecewise discontinuous linear polynomials. Numerical solutions on uniform meshes in space are computed with the constant time step $\tau =0.001$ and  using a $64\times 64$ uniform spatial mesh with DoFs $24576$. 

%%%%%%%%%%%%%%%%%%%%%%%%%%%%%%%%%%%%%%%%%%%%%%%%%%%%%%%%%%%%%%%%%%%%%%%%%%%%%%%
\subsection{Sheering Flow}
\label{ex2}

\begin{figure}[htb!]
\centering
\includegraphics[width=0.32\textwidth]{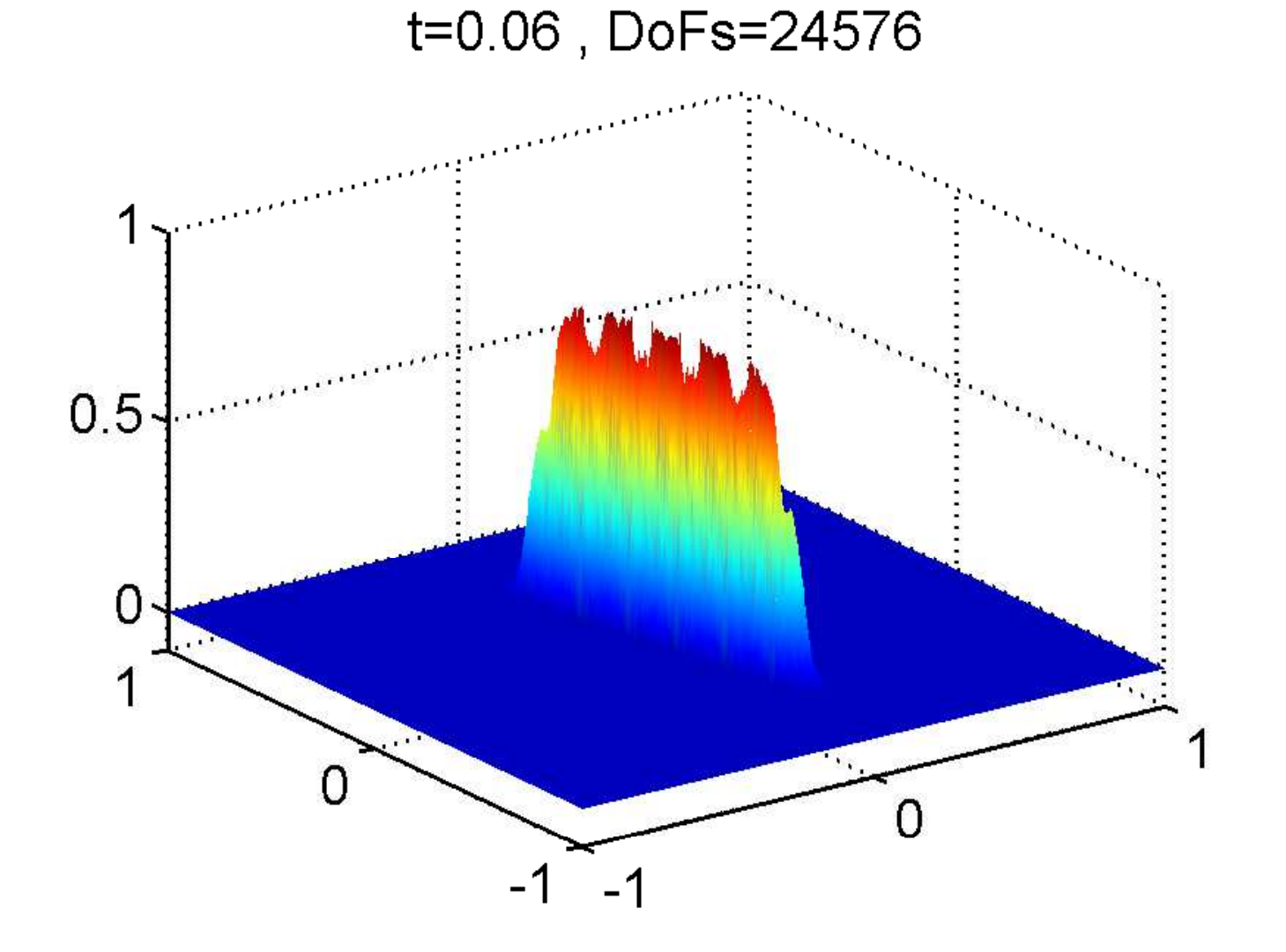}
\includegraphics[width=0.32\textwidth]{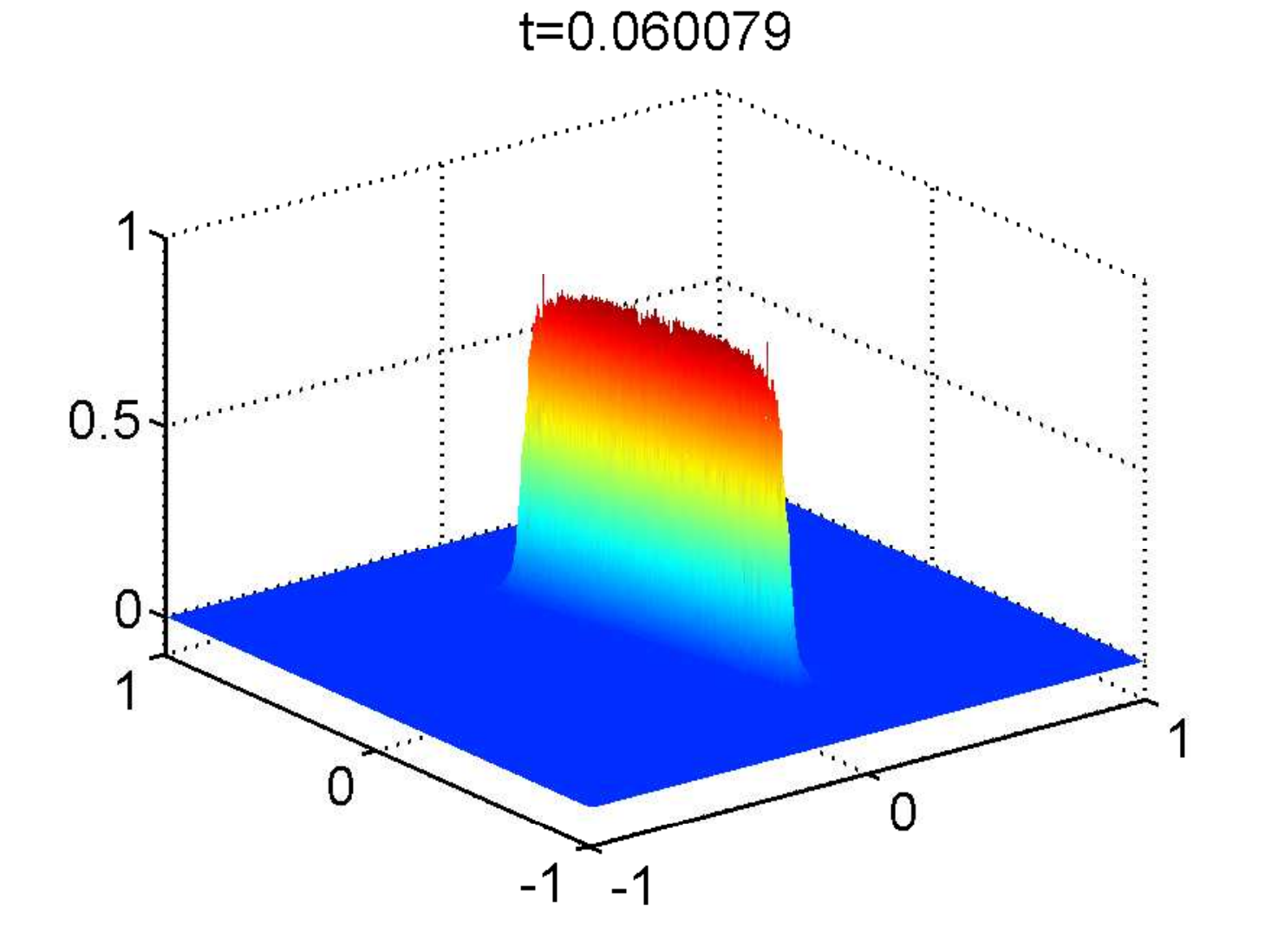}
\includegraphics[width=0.32\textwidth]{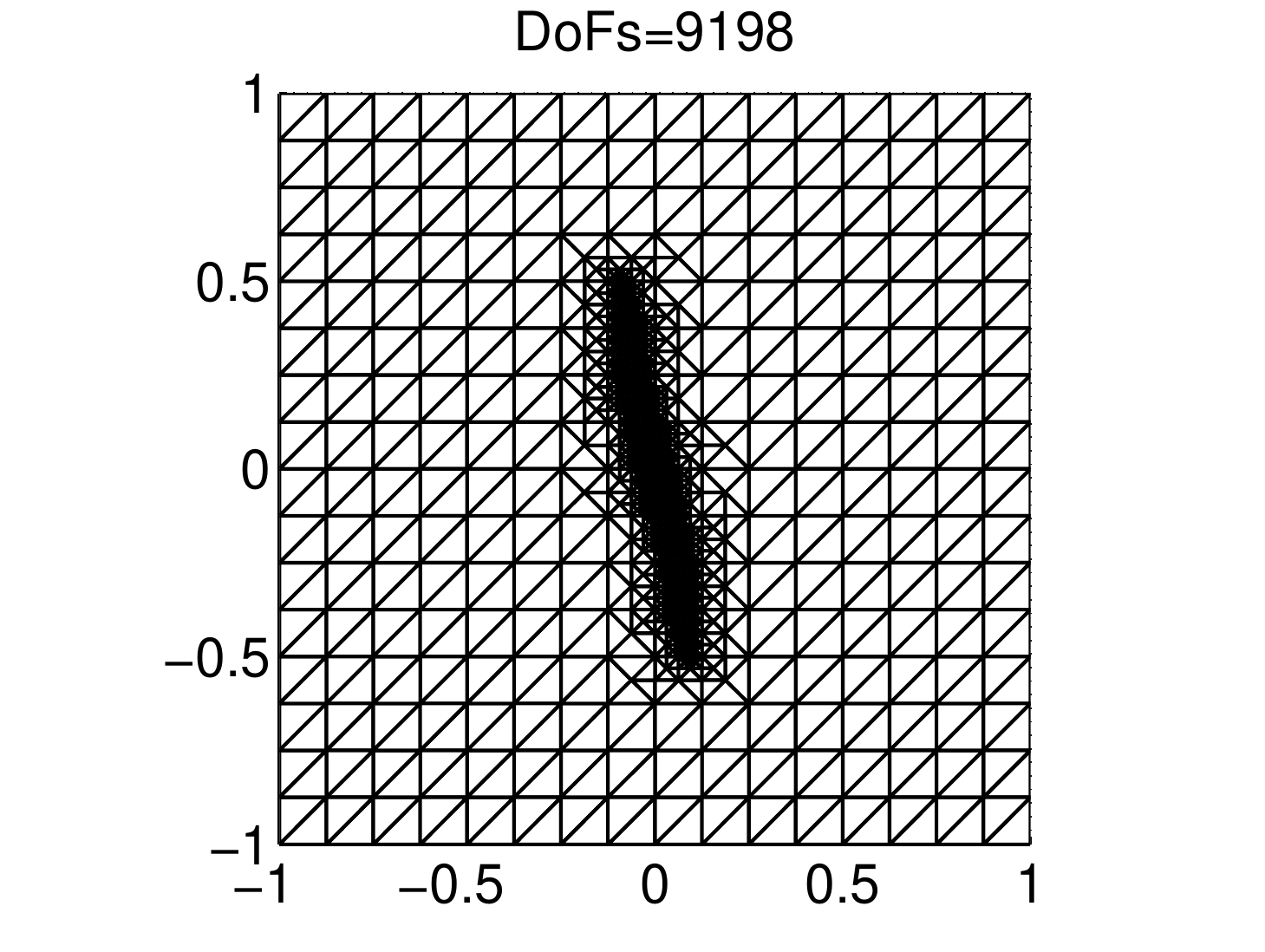}
\caption{Sheering Flow: Solution profiles at final time obtained by uniform (left) and adaptive (middle) schemes, and adaptive mesh at final time (right)\label{plot1}}
\end{figure}

We consider \eqref{advac} with the sheering velocity field $\mathbf{V}=(0,-100 x)$, and with the initial condition as $1$ on $[-0.1,0.1]^2$ otherwise $0$ \cite{liu12}. In Fig.~\ref{plot1}, left, the unphysical oscillations of the solution on uniform mesh can be clearly seen. The oscillations are damped out by the AMOT algorithm, in Fig.~\ref{plot1}, middle, and adaptive mesh is concentrated in the region where the sharp layers occur.

\begin{figure}[htb!]
\centering
\includegraphics[width=0.35\textwidth]{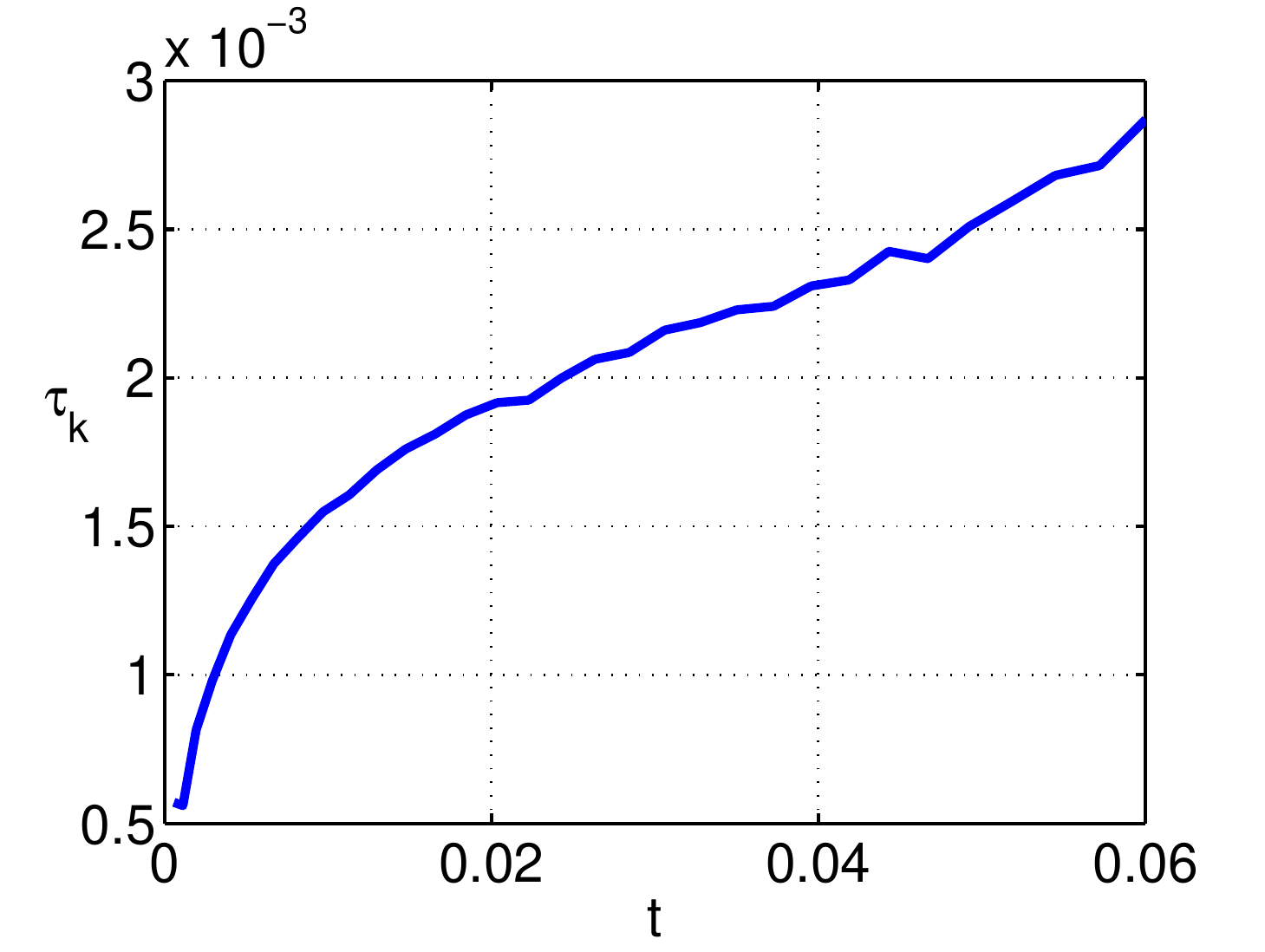}
\includegraphics[width=0.35\textwidth]{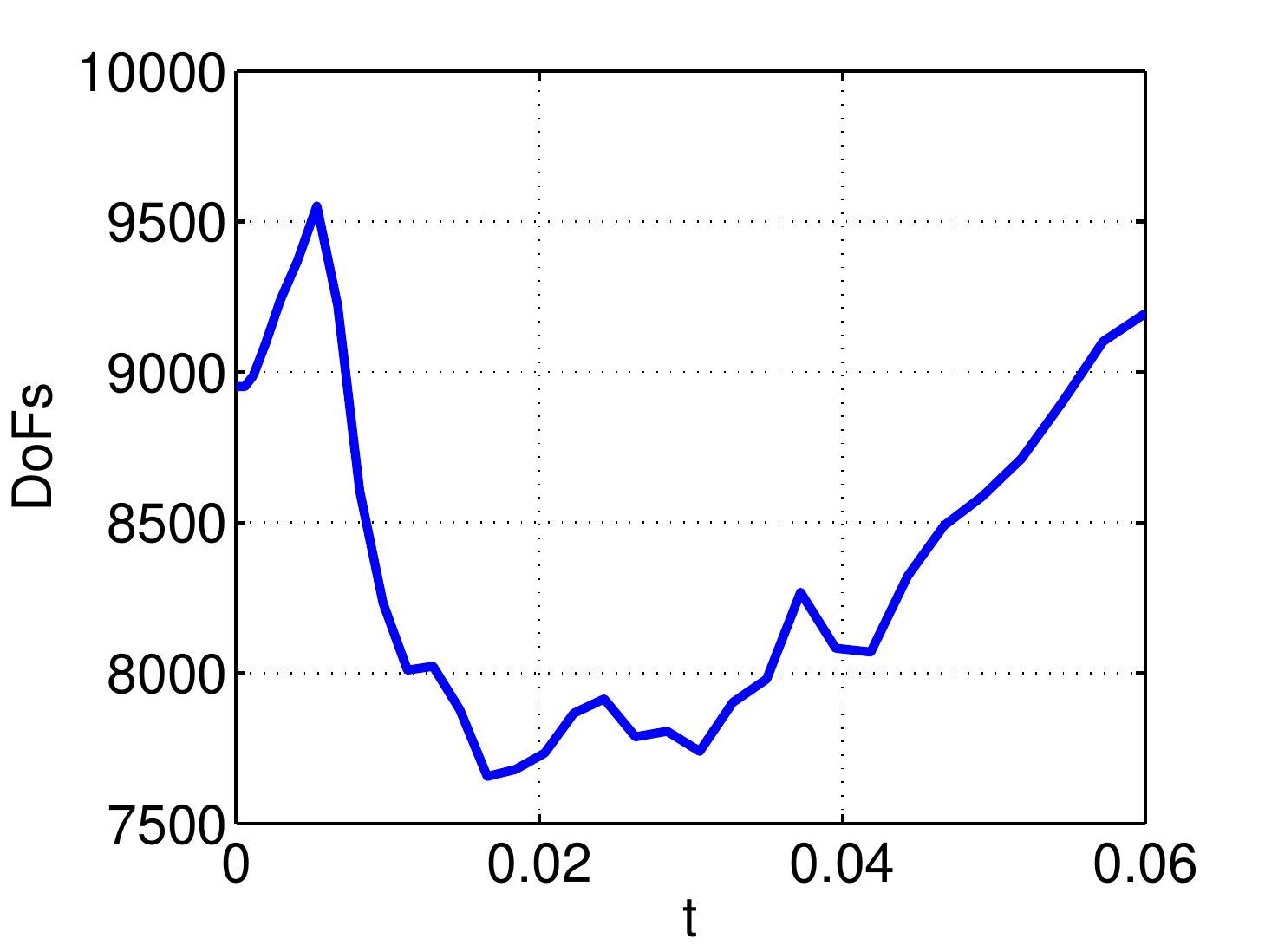}
\caption{Sheering Flow: Evolution of time step size (left) and DoFs (right)  \label{tdt1}}
\end{figure}

The refinement and coarsening of AMOT algorithm works well as shown in Fig.~\ref{tdt1}, right. The mesh becomes finer at the very beginning and then, gets coarser around $t=0.02$ as the size of the interior layer becomes smaller due to the sheering and the time step-size increases monotonically, Fig.~\ref{tdt1}, left.

%%%%%%%%%%%%%%%%%%%%%%%%%%%%%%%%%%%%%%%%%%%%%%%%%%%%%%%%%%%%%%%%%%%%%%%%%%%%%%%
\subsection{Expanding Flow }
\label{ex1}

\begin{figure}[htb!]
\centering
\includegraphics[width=0.32\textwidth]{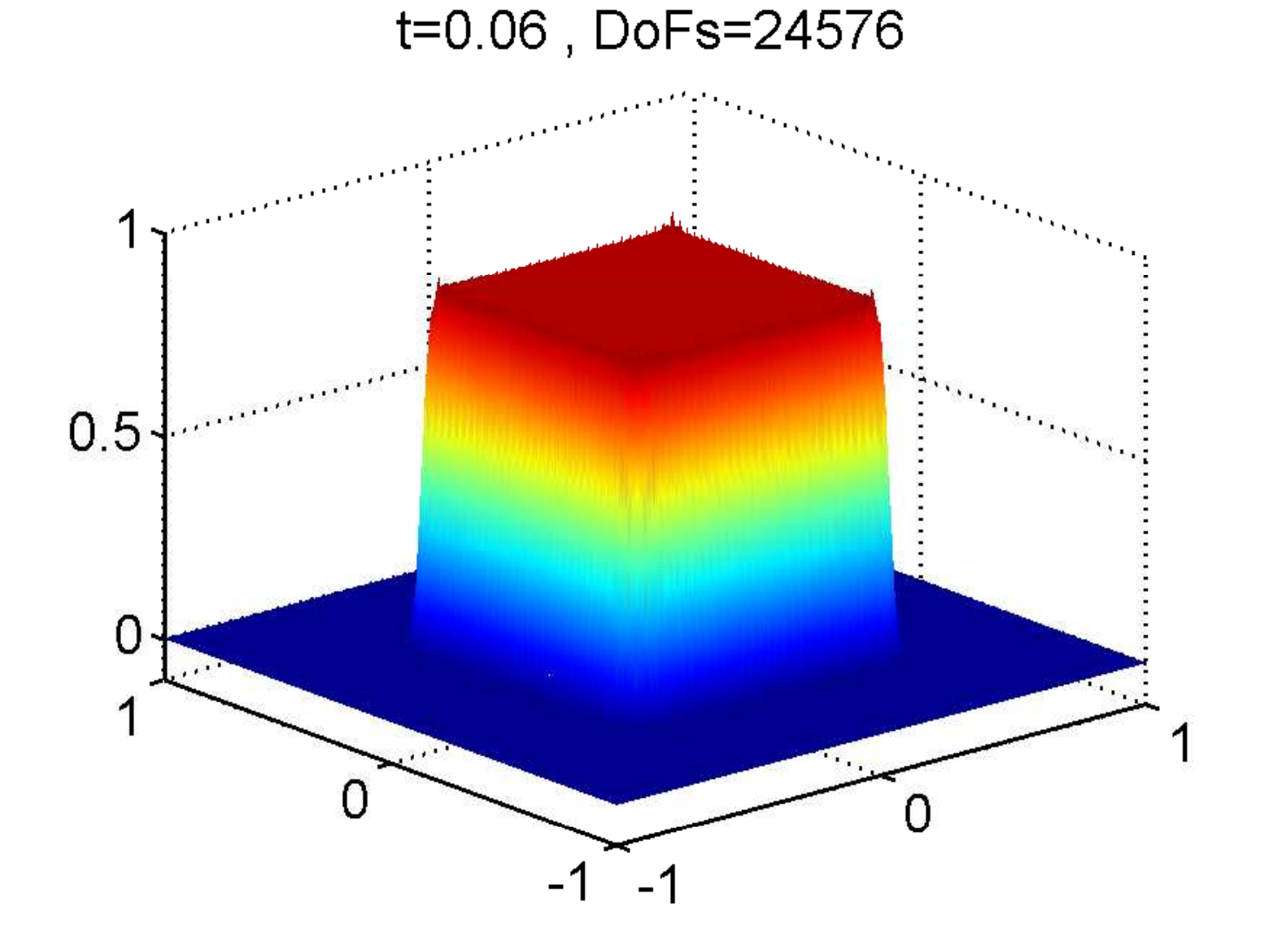}
\includegraphics[width=0.32\textwidth]{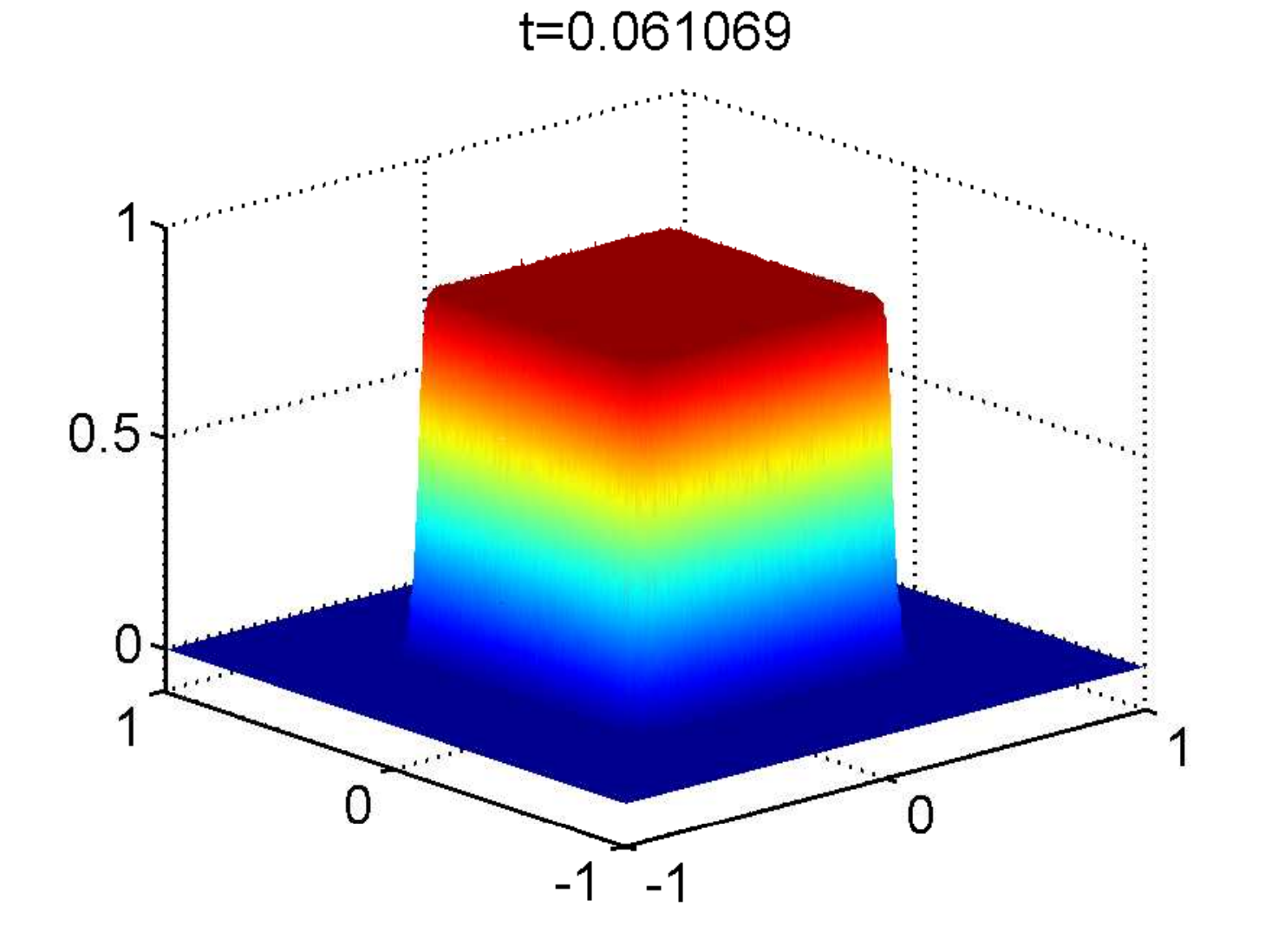}
\includegraphics[width=0.32\textwidth]{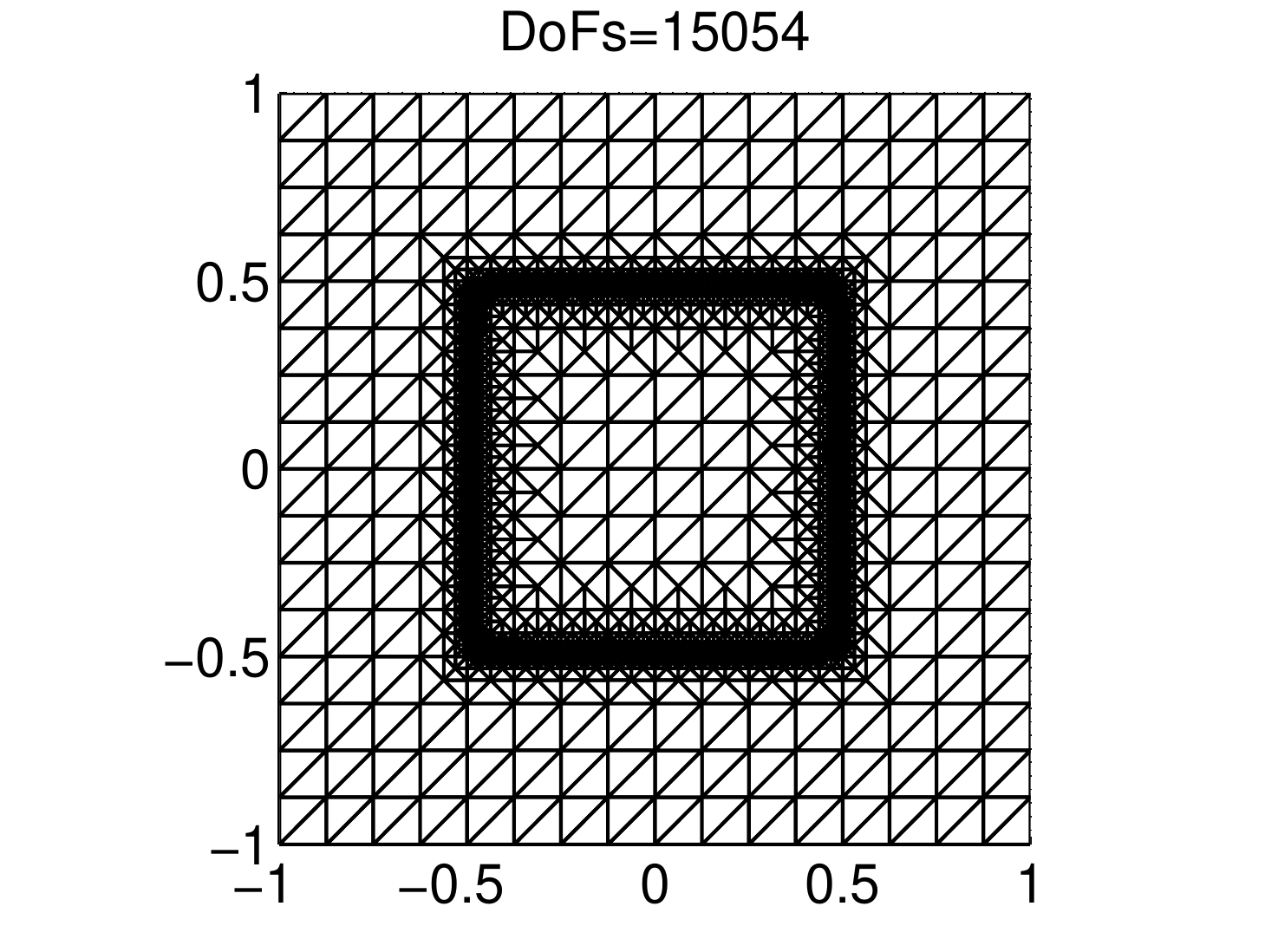}
\caption{Expanding Flow: Solution profiles at final time obtained by uniform (left) and adaptive (middle) schemes, and adaptive mesh at final time ( right)\label{plot2}}
\end{figure}

\begin{figure}[htb!]
\centering
\includegraphics[width=0.35\textwidth]{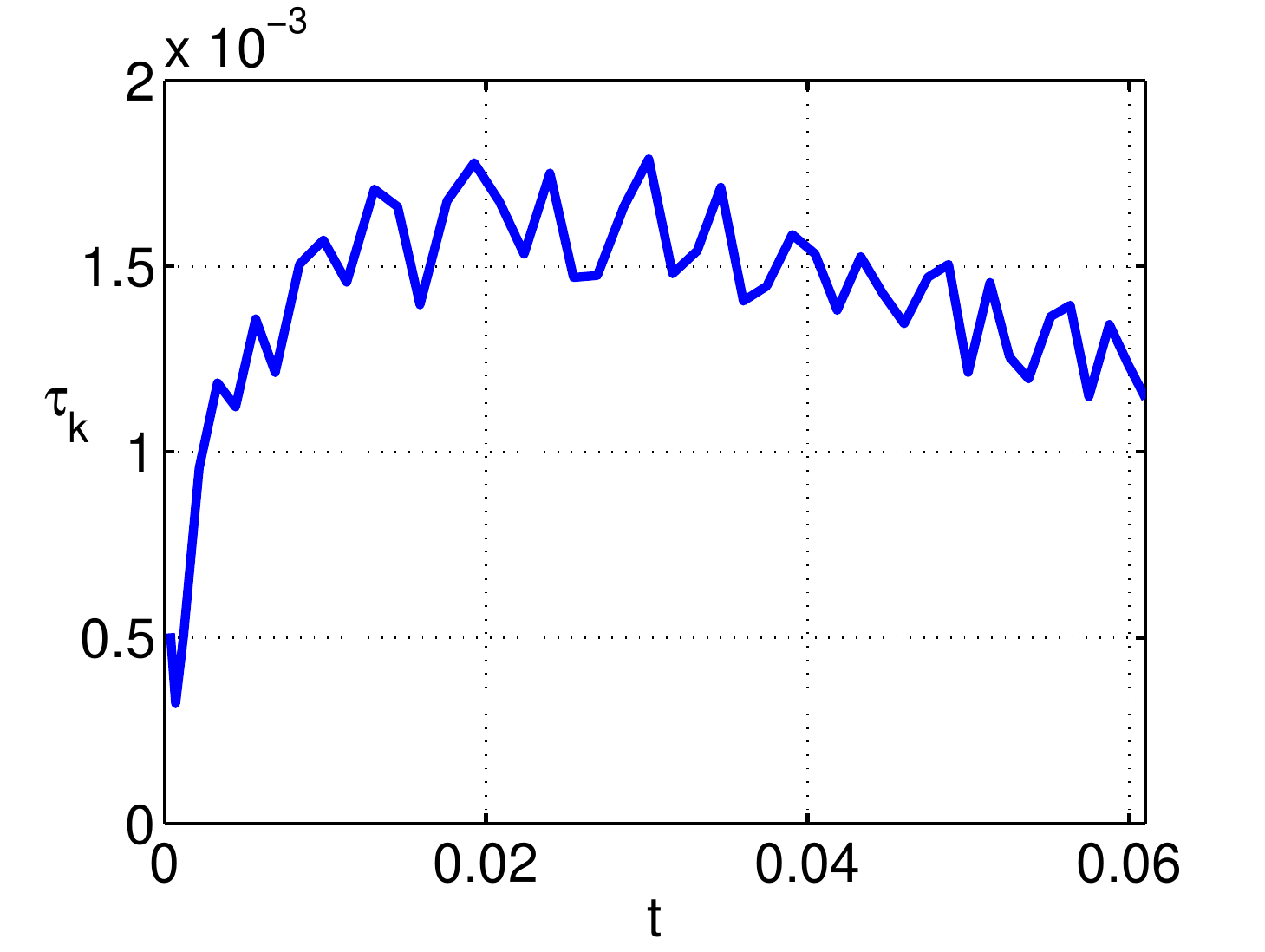}
\includegraphics[width=0.35\textwidth]{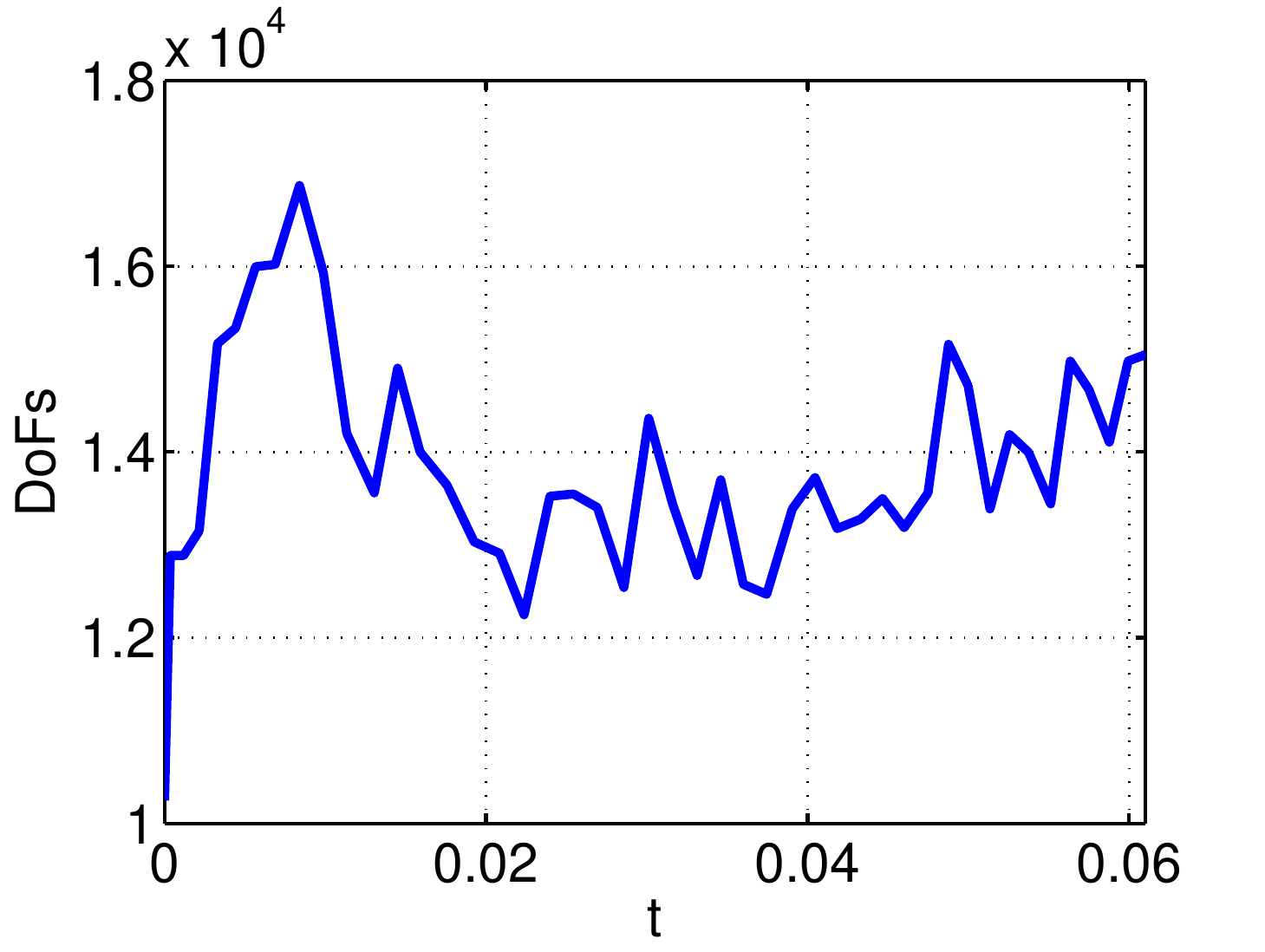}
\caption{Expanding flow: Evolution of time step sizes (left) and DoFs (right) \label{tdt2}}
\end{figure}

As the second example, we consider the expanding velocity field $\mathbf{V}=(10 x,10 y)$. The initial condition is taken as $1$ in the square $[-0.3,0.3]^2$ and $0$ otherwise \cite{liu12}. The unphysical oscillations are damped again in Fig.~\ref{plot2}, middle. The mesh refined slightly and time step-size increases at the beginning, and then refinement and coarsening proceed simultaneously, Fig.~\ref{tdt2}. Time step-size slightly decreases after $t=0.02$ following refinement/coarsening.

%%%%%%%%%%%%%%%%%%%%%%%%%%%%%%%%%%%%%%%%%%%%%%%%%%%%%%%%%%%%%%%%%%%%%%%%%%%
%%%%%%%%%%%%%%%%%%%%%%%%%%%%%%%%%%%%%%%%%%%%%%%%%%%%%%%%%%%%%%%%%%%%%%%%%%%
\section*{Acknowledgments}
This work has been partially supported METU Research Fund  Project BAP-07-05-2013-004 and by Scientific Human Resources Development Program (\"OYP) of the Turkish Higher Education Council (Y\"OK).

%%%%%%%%%%%%%%%%%%%%%%%%%%%%%%%%%%%%%%%%%%%%%%%%%%%%%%%%%%%%%%%%%%%%%%%%%%%

\ifx\undefined\bysame
\newcommand{\bysame}{\leavevmode\hbox to3em{\hrulefill}\,}
\fi

%%%%%%%%%%%%%%%%%%%%%%%%%%%%%%%%%%%%%%%%%%%%%%%%%%%%%%%%%%%%%%%%%%%%%%%%%%%
%%%%%%%%%%%%%%%%%%%%%%%%%%%%%%%%%%%%%%%%%%%%%%%%%%%%%%%%%%%%%%%%%%%%%%%%%%%

\end{document}